\documentstyle[12pt]{article}

\newcommand{\argmax}{\mathop{\rm argmax}\limits}

\newcommand{\const}{\mathop{\rm const}\limits}

\newcommand{\Law}{\mathop{\rm Law}\limits}

\begin{document}

\begin{center}

{\bf PARAMETRIC DENSITY-BASED OPTIMIZATION  } \\

\vspace{4mm}

{\bf OF PARTITION IN CLUSTER ANALYSIS, }\\

\vspace{4mm}

{\bf with applications.} \\

\vspace{4mm}

 $ {\bf E.Ostrovsky^a, \ \ L.Sirota^b,  \ \ A.Zeldin^c.} $ \\

\vspace{4mm}

$ ^a $ Corresponding Author. Department of Mathematics and Computer science, \\
Bar-Ilan University, 84105, Ramat Gan, Israel.\\

\vspace{4mm}

E - mail: \ galo@list.ru \  eugostrovsky@list.ru\\

\vspace{4mm}

$ ^b $  Department of Mathematics and computer science. Bar-Ilan University, 84105,\\ Ramat Gan, Israel.\\

\vspace{4mm}

E - mail: \ sirota3@bezeqint.net\\

\vspace{4mm}

$ ^c $  Research and consulting  officer, the Ministry of Immigrant Absorption, Israel. \\

\vspace{4mm}

E - mail: \ anatolyz@moia.gov.il\\

\vspace{5mm}
                    {\bf Abstract.}\\

 \end{center}

 \vspace{4mm}

  We developed an optimal in the natural sense algorithm of partition in cluster analysis based on the densities of observations
in the different hypotheses. These densities may be characterized, for instance, as the multivariate so-called "quasi-Gaussian distribution".\par
 We describe also the possible applications in technical diagnosis, demography and philology. \par

\vspace{4mm}


\vspace{4mm}


\vspace{4mm}

\section{Introduction. Statement of the problem. Notations. Definitions. }

\vspace{3mm}

   Let us assume $ (\Omega, \cal{B}, {\bf P}) $ a probability space and
 $ (X, \cal{A}, \mu) $ to be a measurable space with sigma-finite measure $ \mu. $\par

  Let $ \{ f_j \} = \{  f_j(x) \}, \ j = 0,1,2,\ldots, N, \ x \in X  $ be a {\it family }
 of probability  densities, i.e. the family of  measurable non - negative functions $ f_j: X \to R $
as

$$
\int_X f_j(x) \ d \mu(x) = 1.
$$
 In the sequel

 $$
   \int  = \int_X, \hspace{6mm} \int f d \mu =  \int_X f(x) \ d \mu(x) = \int_X f(x) \ \mu(dx).
 $$

 Let also $ \xi = \xi(\omega), \ \omega \in \Omega $ be a random variable(r.v.) with values in the set (space) $  X. $
 It may be a random vector or even a random process or fields etc.\par

\vspace{3mm}

{\bf  Definition 1.1.}  We define  the following family of predicates (hypotheses)
 $ H_j, \ j = 0,1,2,\ldots, N: $ the statement $ H_j $ imply that the r.v. $ \xi $ has a density $  f_j(x): $

\vspace{3mm}

$$
H_j: \hspace{4mm} \nu_{\xi}(A) \stackrel{def}{=}  {\bf P} (\xi \in A) = \int_A f_j(x) \ \mu(dx). \eqno(1.1)
$$

\vspace{3mm}

 The condition in (1.1) that the distributions $ \nu_{\xi}(\cdot) $ is absolutely continuous relative to the measure $ \mu(\cdot) $ is
 not essential: arbitrary  finite set of measures may always be dominated. \par

\vspace{3mm}

One of the main problem of the Cluster Analysis (CA)  is the construction of an optimal algorithm in one sense or another
{\it decision rule}, see the classical monographs of M.R.Anderberg  \cite{Anderberg1} and of P.Arabie P., L.J.Hubert L.J., and G. De Soete
 \cite{Arabie1}.\par

\vspace{3mm}

Let us discuss this  more detail. The {\it deterministic } decision rule  $ R = R(G)  $ may be described as
a {\it partition} of a view

$$
 G = \{ G_j \}, \hspace{5mm} G_j \subset X, \hspace{5mm} \cup_{j=0}^N G_j \subset X. \eqno(1.2)
$$

We choose the hypothesis $ H_i $ if and only if $  \xi \in G_i. $ \par

\vspace{3mm}

 This rule is {\it unambiguously} iff

 $$
 \forall (i,k), \ k \ne i \ \Rightarrow G_i \cap G_k = \emptyset, \eqno(1.3)
 $$
and {\it complete}, iff

$$
\cup_{j=1}^N G_j = X. \eqno(1.4)
$$

\vspace{3mm}

 Arbitrary deterministic rule $  R = R(G) $  has errors:

 $$
\alpha_{i,k} = \alpha(R)_{i,k}  \stackrel{def}{=}  \int_{G_i} f_k(x) \ \mu(dx) = \int_{G_i} f_k \ d \mu, \  k \ne i. \eqno(1.5)
 $$
Actually, if the true predicate is $  H_k, $ then $ \alpha_{i,k}  $ is the probability to obtain the hypothesis $  H_i. $\par

 In contradiction, the {\it randomized  } decision rule $ S = S(\phi) $ may be described  as a collection of a  measurable functions
$  \phi = \{  \phi_j \},  \ \phi_j = \phi_j(x), \ x \in X, \ \phi_j: X \to [0,1] $  so that

$$
\alpha_{i,k} = \alpha(S)_{i,k} = \int_X \phi_i(x) \ f_k(x) \ \mu(dx). \eqno(1.6)
$$

 The randomized strategy $  S, $ which includes as a particular case the deterministic rule, is complete iff

$$
\forall x \in X \ \Rightarrow \sum_{j=0}^N \phi_j(x) = 1  \eqno(1.7)
$$
and is unambiguously iff

$$
\forall x \in X, \ \forall k,i: k \ne i \ \Rightarrow \phi_i(x) \phi_k(x) = 0. \eqno(1.8)
$$

\vspace{3mm}

{\it  In what follows we impose on all the considered decisions rules  both the conditions (1.7) and (1.8). }  \par

\vspace{3mm}

 The meaningful sense of the formula (1.6) is evident:  if the true predicate is $  H_k, $
then by means of additional random mechanism independent on $ \Omega $ we admit the hypothesis  $  H_i $ with
probability $ \alpha_{i,k}. $ \par
The probability of a  {\it false alarm }  $ Q_{fa} $ may be expressed through $ \{ \alpha_{i,k} \}: $

$$
Q_{fa} = \sum_{j=1}^N \alpha_{j,0},
$$
as well as the probability of a {\it undetected faults } $  Q_{nd}: $

$$
Q_{nd} = \sum_{j=1}^N \alpha_{0,j}.
$$
 Also the negative sense has probabilities

$$
\overline{\alpha}_{i,i} = 1 - \alpha_{i,i};
$$
it represents the probability to {\it reject} the predicate $  H_j $  under condition that exactly took place.\par
 The statement and solving  of different optimization problems formulated in the CA terms $  \{ \alpha_{i,k}(S) \} $ see in the
classical monographs  \cite{Anderberg1}, \cite{Arabie1}, as well as applications in statistics were described in the books
\cite{Leman1}, \cite{Rao1}. \par

 For the technical applications e.g. in technical diagnosis see \cite{Barzilowich1}, \cite{Barzilowich2},  \cite{Birger1}, \cite{Minakov1}
Here the predicate $  H_0 $ in the technical diagnosis correspondent to the  normal state of the object. \par

\vspace{3mm}

{\bf  The authors are trying in the present paper to highlight some new problems of optimization concerning decision rules, to solve them
and to discuss new applications, especially, in philology. }\par

\vspace{3mm}

 Note that in the statistics the statement of an optimization problem looks as a rule as a minimax one \cite{Leman1}, \cite{Rao1}.
 The case when the domains $ G_j $  have a parallelepipedal form was considered in the article \cite{Minakov1}. This approach is
 traditional in the technical diagnosis, see \cite{Birger1},  and the sizes and the centers of the parallelepipeds are a subject to optimization.\par

The  solution obtained in  \cite{Minakov1} is in general case not complete. \par

\vspace{4mm}

\section{Main result: statement and solving of an optimization problem. }

\vspace{3mm}

{\bf A.  Formation of objective function. } \par

\vspace{3mm}

 Let $  v_{i,k}, \ i,k = 0, 1, \ldots, N, \ k \ne i $  be arbitrary non - negative non - trivial  constants (weights) and
 defined formally as $  v_{i,i} = 0, \ i = 0, 1, \ldots, N. $ \par

 \vspace{3mm}

{\it  We introduce the following objective function (more exactly, functional) }

$$
Z = Z(S) = \sum \sum_{i,k = 0}^N v_{i,k} \ \alpha_{i,k}(S). \eqno(2.1)
$$

\vspace{3mm}

 For instance, the objective function  may look like

 $$
Z = Z(S) = \sum \sum_{i,k = 0}^N  \ \alpha_{i,k}(S)
$$
or

$$
Z = Z(S) = v_1 Q_{nd} + v_2 Q_{fa}
$$
etc. \par

 The weight coefficients may be proportional to the priory probabilities of appearance of the different states $  j  $
or economical damage from faults.     \par

\vspace{3mm}

{\bf B.  Statement of the optimization problem. } \par

\vspace{3mm}

The following  statement of optimization problem seems quite natural. \par
{\it Find the minimum  of the functional} $  Z = Z(S): $

$$
Z = Z(S) = \sum \sum_{i,k = 0}^N v_{i,k} \ \alpha_{i,k}(S) \to \min_S \eqno(2.2)
$$
{\it under conditions}

$$
\phi_k(x) \in [0,1];  \ \sum_{i=0}^N \phi_i(x) = 1; \ \forall k,i: k \ne i \ \Rightarrow \phi_i(x) \phi_k(x) = 0, \eqno(2.3)
$$
(constrained optimization). \par
   The problem (2.2) - (2.3) in the case when the decision rule is deterministic may be reduced as follows. Find the partition
 $ G = \{  G_j \} $   of the set $  X  $ such that

 $$
 Z = Z(R) = \sum \sum_{i,k} v_{i,k}  \int_{G_i} f_k(x) \ \mu(dx) \to \min_G \eqno(2.4)
 $$
under natural conditions

$$
 \forall (i,k), \ k \ne i \ \Rightarrow G_i \cap G_k = \emptyset, \hspace{7mm} \cup_{j=1}^N G_j = X. \eqno(2.5)
 $$

{\bf C. Reducing to the transportation problem.} \par

\vspace{3mm}

 Denote

 $$
 g_i(x) = \sum_{k=0}^N v_{k,i} f_k(x), \eqno(2.6)
 $$
then

$$
Z(S) = \sum_{j=0}^N \int_X \phi_j(x) \ g_j(x) \ \mu(dx). \eqno(2.7)
$$

 In particular,  we can write for the deterministic decision rule

 $$
 Z(R) = \sum_{j=0}^N \int_{G_j} g_j(x) \ \mu(dx). \eqno(2.8)
 $$

Doubtless that the functional $  Z = Z(S) $ is linear over the collection of the functions $  \phi = \{  \phi_i(x) \}. $
The discrete approximation of the functional $  Z $ over the  discrete set $ \{  x_r \} $ of the values $  x, \ x \in X $
looks like

$$
Z(R) \approx  Z_{\Delta}(R) = \sum_r \sum_{j=0}^N \phi_j(x_r) \ g_j(x_r) \ \Delta_r.\eqno(2.9)
$$

 We came to the following optimization problem

$$
\sum_r \sum_{j=0}^N \phi_j(x_r) \ g_j(x_r) \ \Delta_r \to \min_{\phi_j(x_r)} \eqno(2.10)
$$
 under conditions

$$
\forall r \ \Rightarrow \phi_j(x_r) \in [0,1]; \hspace{7mm} \sum_{j=0}^N \phi_j(x_r) \Delta_r = 1, \eqno(2.11)
$$
 or correspondingly

 $$
\forall r \ \Rightarrow \phi_j(x_r) \in (\{0 \}, \ \{1 \} );  \hspace{7mm} \sum_{j=0}^N \phi_j(x_r) \Delta_r = 1. \eqno(2.12)
$$
 The problem (2.10) - (2.11)  belongs to the class of the well - known transportation problems of linear programming. It may be considered as an
approximation  for the source problem  (2.2) - (2.3) and may be used in practice.  \par

\vspace{3mm}

{\bf D.  Solving  of the optimization problem. Main result. } \par

\vspace{3mm}

{\bf  Theorem. } {\it  The optimal decision rule exists, it is unique, deterministic and looks like}

$$
G^{0}_j =  \{ x, \ x \in X, \ g_j(x) = \min_k g_k(x) \}. \eqno(2.13)
$$

 {\it Herewith }

$$
Z( \{  G^0_j  \} ) = \min_{ \{G_j\} } Z( \{  G_j  \} ) = \int_X \min_j g_j(x) \ \mu(dx). \eqno(2.14)
$$

 \vspace{3mm}

 {\bf Proof.} The equality

$$
Z( \{  G^0_j  \} ) =  \int_X \min_j g_j(x) \ \mu(dx)
$$
follows immediately from the definition of the partition $ G^0 = \{  G^0_j  \}.   $ \par

 Let now $ \phi = \{  \phi_j(x) \} $ be other randomized decision rule satisfying the conditions of unambiguousness
and completeness.  We have:

$$
Z( \{\phi_j \} ) = \int \sum_{j=0}^N \phi_j(x) \ g_j(x) \ \mu(dx)  \ge \int \sum_{j=0}^N \phi_j(x) \ \min_k g_k(x) \ \mu(dx) =
$$

$$
\int \ \min_k g_k(x) \ \mu(dx) = Z( \{  G^0_j  \} ), \eqno(2.15)
$$
as long as $ \phi_k(x) \ge 0  $ and $ \sum_j \phi_j(x) = 1. $\par

\vspace{3mm}

 This completes the proof of our theorem.\\

\vspace{3mm}

\section{Quasi - Gaussian distributions. Application in philology. }

\vspace{3mm}

  We assumed above that the densities  $  f_j = f_j(x)  $ are known.  They are for instance approximately Gaussian in the
technical diagnosis, see, e.g. \cite{Birger1}, \cite{Minakov1}.\par
 We will describe in this section the application in philology, in particular, to represent the possible densities
which may appear  therein.\par

 The new so-called quasi-Gaussian distributions which may appear in demography and philology were discussed in the previous
paper of the authors  \cite{Ostrovsky1}.
These distributions were substantiated by means of characterization properties under some natural conditions. \par

\vspace{3mm}

 Let us discuss this in more detail.\par

\vspace{3mm}

 There exist  many characterizations of a two-dimensional, or, more generally, multidimensional Gaussian (normal)
distributions, with independent coordinates.  For example, a characterization by means of independence of linear functionals  or through the
distribution of sums of coordinates, see the classical textbook of W.Feller \cite{Feller1}, p. 77 , p. 498 - 500;  by means of the
properties of conditional   distributions, \cite{Albajar1}, \cite{Kotlarski1}; a characterization by means of the properties of  order statistics
\cite{Jian1}; a characterization by means of some inequalities \cite{Bobkov1}, \cite{Kac1} etc., see also the reference therein.\par
 The famous monograph of A.M.Kagan, Yu.V.Linnik, C.R.Rao \cite{Kagan1} is completely devoted to the characterisation problems
in Mathematical Statistics.\par
 Usually, these characterizations are stable (robust),  \cite{Meshalkin1}, \cite{Zolotarev1}. \par

\vspace{3mm}

 Let us consider the following example.\\

\vspace{3mm}

{\bf Example.} We denote as trivial  for any measurable set $ A, \ A \subset R $ its indicator function by
$ I(A) = I_A(x):   $\par

$$
I_A(x) = 1, \ x \in A; \hspace{5mm} I_A(x) = 0, \ x \notin A.\eqno(3.0)
$$

Let us introduce a {\it family} of  functions

$$
\omega_{\alpha}(x) = \omega_{\alpha}(x; C_1, C_2) := C_1 \ |x|^{\alpha(1)} \ I_{(-\infty,0)}(x) + C_2 \ x^{\alpha(2)} \ I_{ (0,\infty)}(x),
$$

$$
x \in R, \ C_{1,2}= \const \ge 0, \ \alpha = \vec{\alpha}  = (\alpha(1), \alpha(2)),  \ \alpha(1), \alpha(2) = \const > -1, \eqno(3.1)
$$
so that $ \omega_{\alpha}(0) = 0,  $ and a family of a correspondent probability densities of a form

$$
g_{\alpha, \sigma}(x) = g_{\alpha, \sigma}(x; C_1, C_2)  \stackrel{def}{=} \omega_{\alpha}(x; C_1, C_2) \ f_{\sigma}(x). \eqno(3.2)
$$
 Since

$$
I_{\alpha(k)}(\sigma) := \int_0^{\infty} x^{\alpha(k)} \exp \left( -x^2/(2 \sigma^2)  \right) \ dx = 2^{(\alpha(k) - 1)/2  } \
\sigma^{(\alpha(k) + 1)} \ \Gamma((\alpha(k) + 1)/2),
$$
where $  \Gamma(\cdot) $ is ordinary Gamma function, there is the interrelation between the constants $ C_1, C_2:  $

$$
C_1 \ I_{\alpha(1)}(\sigma) + C_2 \ I_{\alpha(2)}(\sigma) = \sigma \ (2 \pi)^{1/2}, \eqno(3.3)
$$
has only one degree of freedom.  In particular, the constant $  C_1 $ may be equal  to zero;
in this case the r.v. $ \xi $ possess  only non - negative values.\par

\vspace{3mm}

{\it  We will denote by $  C_i, K_j $ some finite non - negative constants that are not  necessary to be
the same in different places. }\\

\vspace{3mm}

{\bf  \ Definition 3.1. } The one - dimensional distribution of a r.v. $  \xi   $ with density function of a view
$ x \to g_{\alpha, \sigma}(x - a; C_1, C_2), \ a = \const \in R $ is said  to be quasi - Gaussian or equally quasi - normal. Notation:

$$
\Law(\xi) = QN(a,\alpha,\sigma, C_1, C_2). \eqno(3.4)
$$

\vspace{3mm}
Let us explain the "physical" sense of introduced parameters of these distributions.
 The value $  "a" $ in (3.2) may be called  {\it  quasi - center  } by analogy with normal distribution;  the value $ "\alpha" $
expresses the degree of concentration of this distribution about the center  and the value of
 $ "\sigma" $ which may be  called {\it  quasi - standard } of the r.v. $  \xi  $  expressed alike in the classical Gaussian
r.v.  the degree of scattering.   \\

\vspace{2mm}
 Note that there are some grounds  to accept that the deviation of the point of put-down (landing of air-plane) from the central line of the
landing strip has a quasi-Gaussian distribution, see \cite{Mirzachmedov1},  \cite{Mirzachmedov2}.\par

\vspace{2mm}

  Many properties of these distributions are previously studied in \cite{Ostrovsky1}: moments, bilateral tail behavior etc. In particular,
 it is proved that  if the r.v. $  (\xi, \eta) $ are independent and both have the quasi-Gaussian distribution
with parameters $  a = 0, \ b = 0 $ ("quasi - centered" case):

$$
\Law(\xi) = QN(0,\alpha,\sigma, C_1, C_2), \hspace{5mm} \Law(\eta) = QN(0,\beta,\sigma, C_3, C_4) \eqno(3.5)
$$
may occur with different parameters $ \alpha \ne \beta, \ C_1 \ne C_3, C_2 \ne C_4  $ but with the same value of the standard
$  \sigma, \ \sigma > 0, $ then their polar coordinates  $ (\rho, \zeta)  $ are also independent. \par

  The opposite conclusion was also proved in \cite{Ostrovsky1}: the  characterization of quasi - Gaussian distribution in the demography and
 philology: if the  polar and Decart (cartesian) coordinates  are independent, then under some natural conditions
 the random variables  $ \xi, \eta $  have quasi-Gaussian distribution, and is explained why this property denotes this distribution
 of the words parameters in many languages.  \par

 \vspace{3mm}

 It is possible to  generalize our distributions  on the  multidimensional case. Actually, let us consider the  random vector
$ \xi = \vec{\xi} =  (\xi_1, \xi_2, \ldots, \xi_d) $  with the density

$$
f_{\xi}(x_1, x_2, \ldots, x_d) = G( x_1, x_2, \ldots, x_d; \vec{\alpha}, \vec{\sigma}, \vec{ C_1 },  \vec{C_2  } ) \stackrel{def}{=}
$$

$$
\prod_{j=1}^d  g_{\alpha_j, \sigma_l}(x_j; C_1^{(j)}, C_2^{(j)}), \eqno(3.6)
$$
where $ \alpha_j > -1, \ \sigma_j = \const  > 0, \ C_i^{(j)}  = \const \ge 0, $

$$
C_1^{(j)} \ I_{\alpha^{(j)}(1)}(\sigma_j)  + C_2^{(j)} \ I_{\alpha^{(j)}(2)}(\sigma_j) = \sigma_j \ (2 \pi)^{1/2}, \eqno(3.7)
$$

\vspace{3mm}

  The multidimensional version of our theorem is as follows, \ see \cite{Ostrovsky1}, proposition 3.1:

\vspace{3mm}

{\it Assume that all the standards $ \sigma_j = \sigma $ do not depend on the number $  j. $ Then
the (Cartesian) coordinates of the vector $  \vec{\xi}, $ i.e. the random variables $  \{  \xi_j \} $ are common and independent
and so are their polar coordinates. \par
 The contrary is also true:  if
the Cartesian and polar coordinates of the vector $  \vec{\xi} $ are commonly independent and the random variables $  \{  \xi_j \} $
are regularly distributed, then its density has a form (3.6), with the same standards } $ \sigma. $ \par

\vspace{3mm}

 The knowledge of the densities' form $  f_j = f_j(x) $ of possible distribution $  \xi $ give us a huge advantage
for clusterisation;  but we need to describe the method of parameters measurement (estimation). \par

\vspace{3mm}

\section{Estimation of parameters of quasi-Gaussian distribution.}

\vspace{3mm}

{\bf Definition 4.1  of a weight quasi-Gaussian distributions. } \par

 \vspace{3mm}

 Let $ W_k, \ k = 1,2, \ldots, N $ be positive numbers (weights) such that $ \sum_{k=1}^N W_k = 1. $ We define the {\it weight}
or {\it mixed} quasi - Gaussian distribution  by means of multivariate density like

$$
G^{(W)}( x_1, x_2, \ldots, x_d) = G^{(W)} \left( x_1, x_2, \ldots, x_d; \{ a_j^{(k)} \}, \{ \alpha_j^{(k)} \}, \{ \sigma_j^{(k)} \}, \{ C_1^{(k)} \}\right) \stackrel{def}{=}
$$

$$
\sum_{k=1}^N  W_k \ G \left( x_1-a_1^{(k)}, x_2-a_2^{(k)}, \ldots, x_d-a_d^{(k)}; \vec{\alpha}^{(k)}, \ \{ \sigma_j^{(k)} \}, \ \vec{ C_1 }^{(k)} \right) =
$$

$$
G^{(W)}( \vec{x}, \vec{\theta}), \eqno(4.1)
$$
where

$$
\vec{\theta} = \theta \stackrel{def}{=} \{ N; \{ \vec{a}_d \}, \ \{ \vec{\sigma}_j \}, \  \{ \vec{ C_1 }^{(k)} \}   \}, \ d = \dim X, \ j,k = 0,1,2,\ldots,N.
$$

\vspace{3mm}

A more general form of similar distribution has a discrete component with at the same characterization property:

 $$
 G_0^{(W)}( \vec{x}) := W_0 \delta(\vec{x} - \vec{a_0}) +
 $$

 $$
 \sum_{k=1}^N  W_k \ G \left( x_1-a_1^{(k)}, x_2-a_2^{(k)}, \ldots, x_d-a_d^{(k)}; \vec{\alpha}^{(k)}, \vec{\sigma}^{(k)}, \ \vec{ C_1 }^{(k)} \right), \eqno(4.2)
 $$

$$
W_0, W_1, \ldots, W_N > 0, \ \sum_{k=0}^N W_K = 1,
$$

$ \delta(\vec{x}) $ is the classical Dirac delta function; so that

$$
{\bf P} (\vec{\xi} =  \vec{a_0}) = W_0 > 0.
$$

\vspace{3mm}

 {\it Statement of problem:} given a sample $   \{  \eta_m \}, \ m = 1,2,\ldots, n; $ where $  n >> 1 $  from the weight multivariate
quasi-Gaussian distribution; we need to estimate its parameters: the number of clusters $ N, $ the centers $ \vec{a}_j, $ degrees of
concentrations $ \vec{\alpha} $ etc. \par

\vspace{3mm}

Let us imagine it by means of the {\it demography} analogy.  Here the weights $ W_k $ are proportional to the
 share of $ k^{th} $ city in the general population of some country.\par
  In contradiction, in the {\it philology} the parameters $ \{W_k \} $  are possibly unknown and are
subject to evaluation on a sample.

\vspace{3mm}

Regarding the applications of the developed method in linguistics, let's consider the bunch of words of similar meaning (e.g.
hand, arm, palm, elbow, thumb, finger, to take, to give, to get, to bring, to catch, to hold etc,), so-called "semantic field". These words
are grouped around  a semantic nucleus (here - the notion of hand/arm) and will be considered as a cluster. It may be compared with
other clusters in order to calculate lexical/semantic affinity on the base of the proposed quasi-Gaussian distribution. The results may
suggest the common origin, provided the etymological analysis permits it.

 \vspace{3mm}

 The very same equation (more precisely, system of equations) of maximal likelihood for the parameters estimation has a classical
form:

$$
\hat{\theta}_n = \argmax_{\theta} \sum_{m=1}^n \log  G^{(W)}( \vec{\eta}, \vec{\theta}). \eqno(4.3)
$$

 It is well-known  that the rate of convergence $ \hat{\theta}_n $ to the true value $ \theta_0 $ is $ n^{-1/2}. $
The non-asymptotic deviation

$$
{\bf P}_{\theta,n}(u) := {\bf P} \left( \sqrt{n} || \hat{\theta}_n  - \theta_0 || > u \right)
$$
as

$$
{\bf P}_{\theta,n}(u) \le \exp \left(- K(\theta) \ u^{ \gamma } \right)  \approx \exp \left(- K(\theta_n) \ u^{ \gamma } \right), \ u \ge 1,
\ \gamma = \const > 0 \eqno(4.4)
$$
is studied in \cite{Ostrovsky2}.\par

 Moreover,

 $$
 {\bf P} \left( \hat{N}_n \ne N  \right) \le C_5(\vec{\theta}) \ q^n( \vec{\theta}),  \hspace{6mm}  0 < q^n( \vec{\theta}) < 1. \eqno(4.5)
 $$

  The quasi - centers $ \{  a_j^{(k)} \} $ may be interpreted as coordinates of fundamental human notions:  food, policy, medicine, economic etc. \par

 One of the important advantage of approach offered above is the {\it automatical} measurement of cluster's number $ N \approx \hat{N}_n, $ in
contradiction to the classical methods of cluster analysis, see  \cite{Anderberg1}, \cite{Arabie1}. \par

 Notice that wherein the speed of convergence  $ \hat{N}_n  $ to the true value of number clusters $  N  $ is very   hight, see (4.5). \par

 We emphasise also that we do not used arbitrary distance between the values $ \eta_m. $ \par

 \vspace{3mm}

 This was made possible only because we deduced the possible form of the distributions $  f_j(x) $ in the parametric form.\par

\vspace{3mm}

 The classification based on the  mixed quasi-Gaussian distribution may be useful, for example, in learning a foreign language.\par

\vspace{3mm}

Needless to say, this approach requires an experimental verification.\par

\end{document}